\documentclass[10pt]{article} 

\usepackage{multirow}
\usepackage{amsmath}
\usepackage{amsopn}
\usepackage{amsthm}
\usepackage{amsfonts}
\usepackage[T1]{fontenc}
\usepackage[latin1]{inputenc}
\usepackage{graphicx,color}
\usepackage{epsfig}
\DeclareGraphicsExtensions{.eps,.ps}
\usepackage{amssymb}
\usepackage[frenchb,english]{babel}
\usepackage{url}
\usepackage[mathscr]{eucal}
\usepackage{amsbsy}
\usepackage{fancyhdr}
\usepackage{hyperref}

\numberwithin{equation}{section}

\newcounter{rq}
\newcommand\remarque{\refstepcounter{rq}\par\rm{\bf Remark \arabic{rq}.}\quad}
\theoremstyle{plain}
\newtheorem*{thma}{Theorem A}

\newcounter{ex}
\newcommand\exemple{\refstepcounter{ex}\par\rm{\bf Example \arabic{ex}.}\quad}

\theoremstyle{plain}
\newtheorem{thm}{Theorem}[section]
\newtheorem{definition}[thm]{Definition}

\newtheorem{proposition}[thm]{Proposition}
\newtheorem{corollaire}[thm]{Corollary}
\newtheorem{lemme}[thm]{Lemma}
\def\one{\hbox{1\hskip -3pt I}}
\def\rd{\mathbb{R}^d}
\def\sd{S^{d-1}}
\def\qed{\hspace*{2mm} \hfill $\Box$\bigskip}
\def\proof{\par\rm{\it Proof.}\quad}

\def\reip{{\bf P}}

\def\reie{{\bf E}}

\def\cvl1{\stackrel{\mathcal{L}_1}{\longrightarrow}}
\def\cvps{\xrightarrow[n\rightarrow\infty]{a.s.}}

\def\cvpsn{\xrightarrow[N\rightarrow\infty]{a.s.}}

\def\cvas{\xrightarrow[]{a.s.}}

\def\cvdn{\xrightarrow[N\rightarrow\infty]{D}}
\def\cvs{\xrightarrow[n\rightarrow\infty]{}}

\def\egenloi{\stackrel{\mathcal{L}}{=}}

\title{Estimation of the spectral measure of multivariate regularly varying
distributions}
\author{
\renewcommand{\thefootnote}{\alph{footnote}}
SHUYAN LIU\,\footnotemark[1]{}
}
\date{}
\begin{document}
\thispagestyle{empty}
\maketitle
\renewcommand{\thefootnote}{\alph{footnote}\,}
\footnotetext[1]{Institut de statistique,
Universit\'e Catholique de Louvain,
Voie du Roman Pays, 20,
B-1348 Louvain-la-Neuve, Belgium.}

\renewcommand{\thefootnote}{\arabic{footnote}}

\begin{abstract}
In the paper, the estimator for the spectral measure of multivariate stable
distributions introduced by 
Davydov and co-workers are extended to the regularly varying
distributions. The sampling method is modified to optimize the rate of
convergence of estimator. An estimator of the total mass of spectral
measure is proposed. The consistency and the asymptotic normality of estimators
are proved.
\end{abstract}

\selectlanguage{english}
\medskip 
\noindent\textit{AMS Classifications : } 60B99, 60E07, 62F10, 62F12.

\smallskip\noindent \textit{Key words and phrases : } regularly varying
distributions, estimation of parameters, stable distributions, spectral
measure

\newpage


\section{Introduction}
A random $\rd$-valued vector $X$ has a regularly varying distribution with
characteristic exponent $\alpha>0$ if there exists a finite measure $\sigma$ in
the unit sphere $\sd=\left\{x~\left|~\|x\|=1, x\in\rd\right.\right\}$ such that
$\forall B\in\mathcal{B}(\sd)$ with $\sigma(\partial B)=0$,
\begin{equation}\label{vrintro}
\lim\limits_{x\rightarrow\infty}\frac{x^\alpha}{L(x)}\reip\left\{\frac{X}{\|X\|}
\in B, \|X\|>x\right\}=\sigma(B),
\end{equation}
where $L$ is a slowly varying function, i.e., $\frac{L(\lambda
x)}{L(x)}\rightarrow 1$ as $x\rightarrow\infty, \; \forall \lambda >0$. Here the
notation $\|\cdot\|$ denotes Euclidean norm. The
measure $\sigma$ is called {\em spectral measure}, and $\alpha$ is called simply
{\em tail index}. The unit sphere $\sd$ will be simply denoted by $S$. The fact
that $X$ has a regularly varying distribution
with tail index $\alpha$ and spectral measure $\sigma$ will be noted later by
"$X\in \mbox{RV}(\alpha,\sigma)$". Regular variation condition appears
frequently in the studies of the limit theorem for normalized sums of i.i.d.
random terms, see e.g. \cite{Rvaceva62} and \cite{Meerschaert01}, and the
extreme value theory, see e.g. \cite{Resnick87}. Regular variation is necessary
and sufficient conditions for a random $\rd$-valued vector belongs to the domain
of attraction of a strictly $\alpha$-stable distribution, if $\alpha\in(0,2)$,
see e.g. \cite{Araujo80}.

There are various characterizations of the property $X\in
\mbox{RV}(\alpha,\sigma)$ (see e.g. \cite{Mikosch03}). We give here an
equivalent definition.
\begin{definition}\label{def2}
A random $\rd$-valued vector $X\in \mbox{RV}(\alpha,\sigma)$ if there
exists a slowly varying function $\tilde L$ such that for all $r>0$ and
$B\in\mathcal{B}(S)$ with $\sigma(\partial B)=0$
\begin{equation}\label{regulier2}
\lim_{n\rightarrow\infty}n\reip\left\{\frac{X}{\|X\|}\in B,
\|X\|>rb_n\right\}=\sigma(B)r^{-\alpha}, 
\end{equation}
where $b_{n}=n^{1/\alpha}\tilde L(n)$.
\end{definition}
It is well known that in $\rd$ the convergence (\ref{regulier2}) are equivalent
to the convergence in distribution of binomial point processes
$\beta_n=\sum^{n}_{k=1}\delta_{X_k/b_n}$ to a Poisson point process
$\pi_{\alpha,\sigma}$ whose intensity measure has a particular form
\cite{Resnick87}. This result is generalized to random elements in an abstract
cone \cite{Davydov08}. Moreover, in a cone which possesses the sub-invariant
norm, the convergence (\ref{regulier2}) with $\alpha\in (0, 1)$ implies that $X$
belongs to the domain of attraction of a strictly $\alpha$-stable distribution
(see Th. 4.7 \cite{Davydov08}).

We are interested in the problem of estimation of the tail index
$\alpha$ and the spectral measure $\sigma$ of a regularly varying distribution.
By using the relation between the stable distributions and the point processes,
Davydov and co-workers (see \cite{Davydov99}, \cite{Davydov00a} and
\cite{Paulauskas03}) proposed a method to estimate the tail index $\alpha$ and
the normalized spectral measure $\sigma$ of the stable distributions in $\rd$.
The objective of this work is to extend this method to the multivariate
regularly varying distribution.

Suppose that we have a sample $\xi_1, \xi_2,\ldots, \xi_N$, taken from a
regularly varying distribution in $\rd$ with unknown tail index $\alpha$ and
unknown spectral measure $\sigma$. We divide the sample into $n$ groups
$G_{m,1},\ldots,G_{m,n}$, each group containing $m$ random vectors. In
practice, we choose 
\begin{equation}\label{setting}
n=[N^r], r\in (0,1), \;\; \mbox{and then}\;\; m=[N/n],
\end{equation}
where $[a]$ stands for the integer part of a number $a>0$. As $N$ tends to
infinity, we have
$nm\sim N$. Let
\begin{equation}\label{defm1}
M^{(1)}_{m,i}=\max\{\|\xi\| \; | \; \xi\in G_{m,i}\}, \ i=1,\ldots,n,
\end{equation}
that is, $M^{(1)}_{m,i}$ denote the largest norm in the group $G_{m,i}$. Let
$\xi_{m,i}=\xi_j=\xi_{j(m,i)}$ where the index $j(m,i)$ is such that 
\begin{equation}\label{defximi}
\|\xi_{j(m,i)}\|=M^{(1)}_{m,i}.
\end{equation}
We set
\begin{equation}\label{defm2}
M^{(2)}_{m,i}=\max\{\|\xi\| \; | \; \xi\in G_{m,i}\backslash\{\xi_{m,i}\}\}, \
i=1,\ldots,n,
\end{equation}
that is, $M^{(2)}_{m,i}$ denote the second largest norm in the same group.
Let us denote
\begin{equation*}
\varkappa_{m,i}=\frac{M_{m,i}^{(2)}}{M_{m,i}^{(1)}}, \;\;\; \
S_n=\sum_{i=1}^n\varkappa_{m,i}
\end{equation*}
and
\begin{equation*}
\hat\alpha_N=\frac{S_n}{n-S_n}.
\end{equation*}

The regular variation condition (\ref{vrintro}) implies
\begin{equation}\label{cdqlrd}
\reip\{\|\xi\|>x\}=x^{-\alpha}\sigma(S)L(x)+o(x^{-\alpha}L(x))
\;\;\;\mbox{as}\;\;\;
x\rightarrow\infty.
\end{equation}
In the following we will need
the stronger relation : for sufficiently large $x$ and for some $\beta>\alpha$
\begin{equation}\label{cdsdord}
\reip\{\|\xi\|>x\}=C_1x^{-\alpha}+C_2x^{-\beta}+o(x^{-\beta}).
\end{equation}

Under the regular variation assumption and the
second-order asymptotic relation (\ref{cdsdord}), the consistency and the
asymptotic normality of the estimator $\hat\alpha_N$ were proved firstly for
$n=m=[\sqrt{N}]$ in \cite{Davydov00a} and then for more general setting
(\ref{setting}) in \cite{Paulauskas03}. We resume these results in the following
theorem.

\begin{thma}\label{thmpaul}(\cite{Paulauskas03})
Let $\xi, \xi_1,\ldots,\xi_N$ be i.i.d. random $\rd$-valued vectors with a
distribution satisfying (\ref{cdqlrd}) and let the numbers $n$ and $m$ satisfy
the relation (\ref{setting}), then
\begin{equation*}
\frac{1}{n}S_n\cvpsn\frac{\alpha}{1+\alpha}.
\end{equation*}
If the distribution of $\xi$ satisfies (\ref{cdsdord}) with
$0<\alpha<\beta\leq\infty$ and we choose
\begin{equation*}
n=N^{2\zeta/(1+2\zeta)-\varepsilon}, \ m=N^{1/(1+2\zeta)+\varepsilon},
\end{equation*}
where $\zeta=(\beta-\alpha)/\alpha$ and $\varepsilon\rightarrow 0$ as
$N\rightarrow\infty$, then as $N\rightarrow\infty$
\begin{equation*}
\frac{\displaystyle\sqrt{n}\left(\frac{1}{n}S_n-\frac{\alpha}{\alpha+1}\right)}{
\displaystyle\left(\frac{1}{n}\sum_{i=1}^n\varkappa_{m,i}^2-\left(\frac{1}{n}
S_n\right)^2\right)^{1/2}}\Rightarrow
\mathcal{N}(0,1).
\end{equation*}
\end{thma}

\vspace{0.5cm}
This paper focus on the estimation of the spectral measure. In \cite{Davydov00a}
an estimator of normalized spectral measure
$\tilde\sigma(\cdot)=\sigma(\cdot)/\sigma(S)$ was proposed as follows. We set 
\begin{equation}\label{deftheta}
\theta_{m,i}=\frac{\xi_{m,i}}{\|\xi_{m,i}\|}, \ i=1,\ldots,n
\end{equation}
where $\xi_{m,i}$ is defined by (\ref{defximi}). Let us denote
\begin{equation}\label{defestms}
\hat{\sigma}_N(\cdot)=\frac{1}{n}\sum_{i=1}^{n}\delta_{\theta_{m,i}}(\cdot).
\end{equation}
Random vectors $\theta_{m,1},\ldots,\theta_{m,n}$ are i.i.d. and it is proved in
\cite{Davydov00a} that $\hat{\sigma}_N(\cdot)$ is consistent considering a fixed
set $B$, that is, $\forall B\in\mathcal{B}(S)$ with $\sigma(\partial B)=0$,
\[\hat{\sigma}_N(B)\cvas\tilde\sigma(B).\]
The asymptotic normality for $\hat{\sigma}_N(B)$ was proved in \cite{Davydov99}.
All these relations were obtained under the assumption that $n=m$. Inspired by
the work in \cite{Paulauskas03} we modify the sampling method of regrouping and
discuss the convergence rate of the estimator of spectral measure for the
general setting (\ref{setting}). By finding a countable collection of the
$\sigma$-continuity sets which is
closed under the operation of finite intersection, we obtain
$\hat{\sigma}_N\stackrel{a.s.}{\Rightarrow}\tilde\sigma \ \mbox{as} \
N\rightarrow\infty$, where $\Rightarrow$ indicate convergence in distribution.

Since $\hat{\sigma}_N$ in (\ref{defestms}) gives the normalized spectral
measure, it remains for us to estimate the total mass $\sigma(S)$. 
Note that the condition (\ref{vrintro}) implies
\[\lim_{x\rightarrow\infty}x^{\alpha}L(x)^{-1}\reip\{\|\xi\|>x\}=\sigma(S)\]
and therefore the value of $\sigma(S)$ depends on the choice of slowly varying
function $L(x)$. Estimation of this function is discussed in \cite{Resnick06}
and \cite{Resnick97}.
Here we assume that the random vector $\xi$ satisfies the condition
(\ref{vrintro}) with $L(x)=1$. That means if $\alpha\in(0,2)$, the law of $\xi$
belongs to the 
normal domain of attraction of an $\alpha$-stable distribution. Let us denote
\begin{equation}\label{defqmi}
q_{m,i}=\frac{M_{m,i}^{(1)}}{m^{1/\alpha}}
\end{equation} 
where $M_{m,i}^{(1)}$ is defined by (\ref{defm1}). The proposed estimator is
defined by
\begin{equation}\label{etrmas}
\widehat{\sigma(S)}_N=\left(\frac{1}{n\Gamma(1-\frac{t}{\alpha})}\sum_{i=1}^n
q_{m,i}^t\right)^{\frac{\alpha}{t}}, \ t>0.
\end{equation}

\begin{figure}[!htbp]
\begin{center} 
\begin{tabular}{cc}
\hspace{-3mm}\includegraphics[width=60mm,
height=40mm]{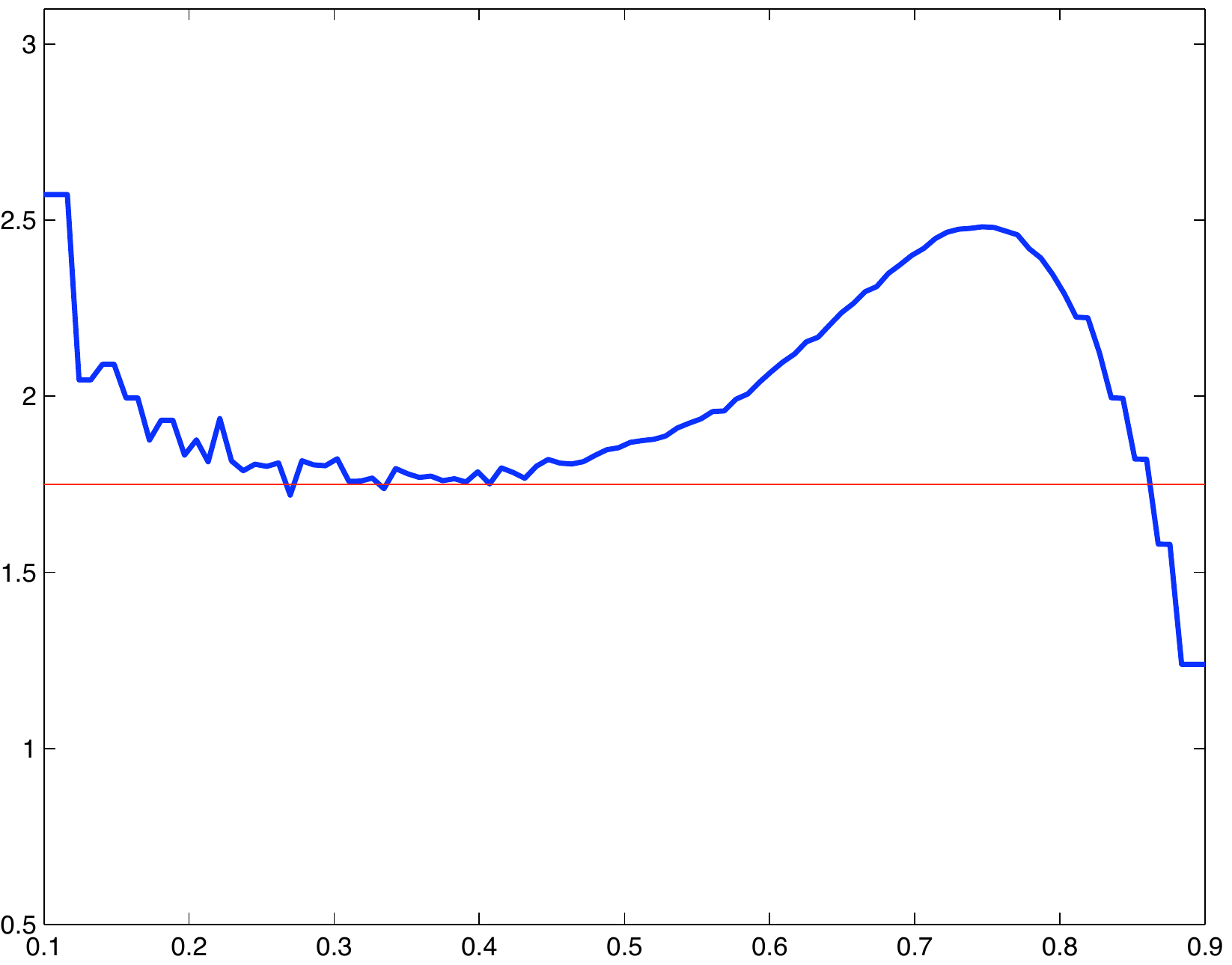}&
\includegraphics[width=60mm, height=40mm]{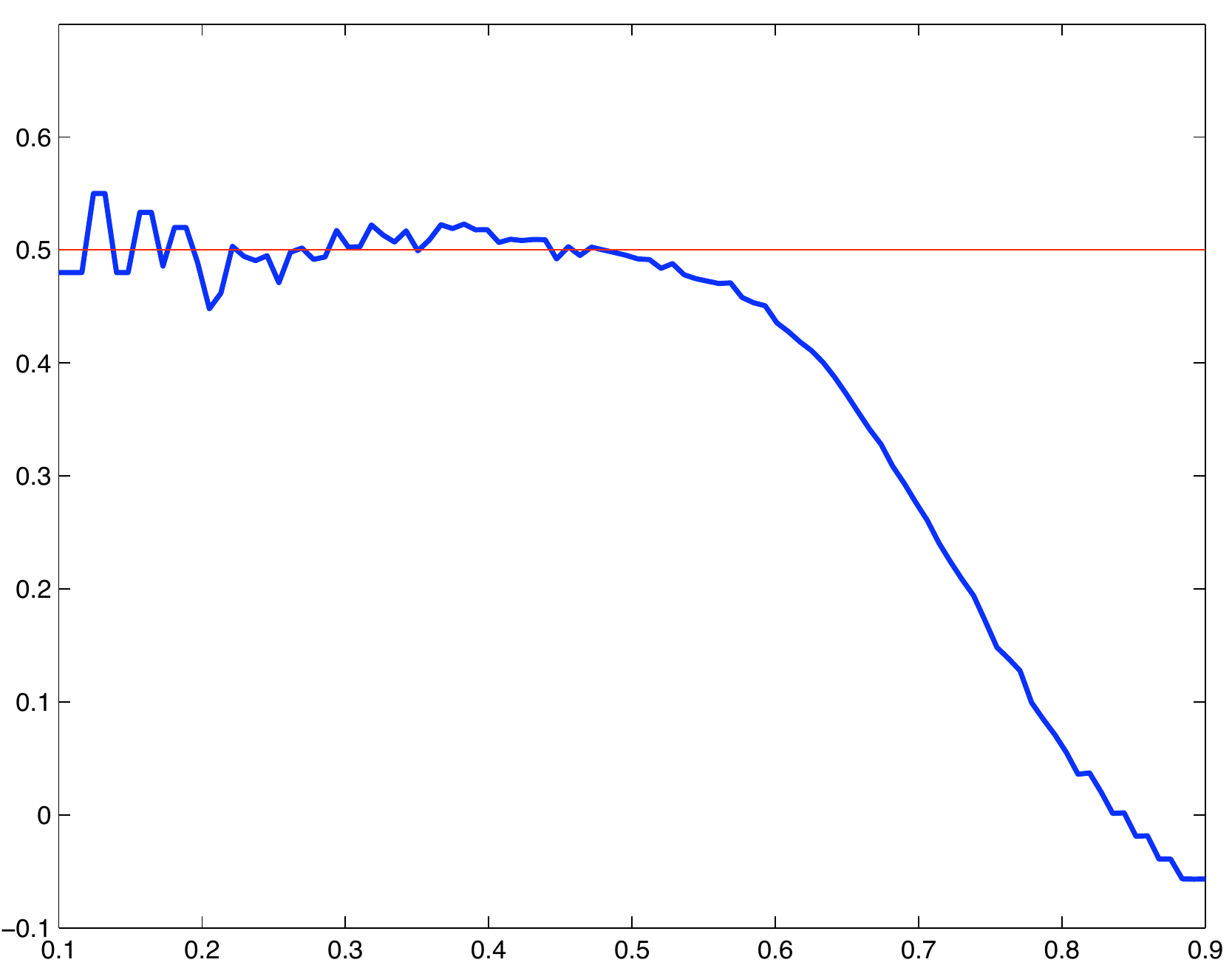}
\end{tabular}
\end{center}
\caption{Estimation results of parameters $\alpha$ (left) and $\rho$ (droite) in
terms of $r$. 
The $x$-axis represent $1-r$. Horizontal line represent the true value.}
\label{influence}
\end{figure}

\exemple
We generated samples from univariate stable distribution with $\alpha=1.75$,
$\sigma(S)=1$ and
$\rho=\frac{\sigma(1)-\sigma(-1)}{\sigma(S)}=0.5$. The sample size is
$100,000$. 
We calculated the estimators $\hat p=\hat p_r$
as a function of $r$. This procedure was repeated $50$ times on the independent
sets of samples. 
Then we plotted $\{(1-r,\bar p_r), 0<r<1\}$ where $\bar p_r$ is the mean of $50$
estimated values
$\hat p_r$. Since the estimator of total mass depends on $\alpha$, we present
here the simulated 
results of the parameters $\alpha$ and $\rho$ in Figure \ref{influence}. The
horizontal line corresponds
to the true value of the parameter $p$. It seems that both plots have the
optimal value of $1-r$ which 
is around $0.4$. In fact by Theorem A and the fact that $\beta=2\alpha$ for a
stable 
random variable, the asymptotically optimal value of $1-r$ is approximately
$1/3+\varepsilon$. 
A similar result for the estimator $\hat\sigma$ is given in Theorem
\ref{normademc}.

\begin{figure}[!htbp]
\begin{center} 
\includegraphics[width=60mm, height=50mm]{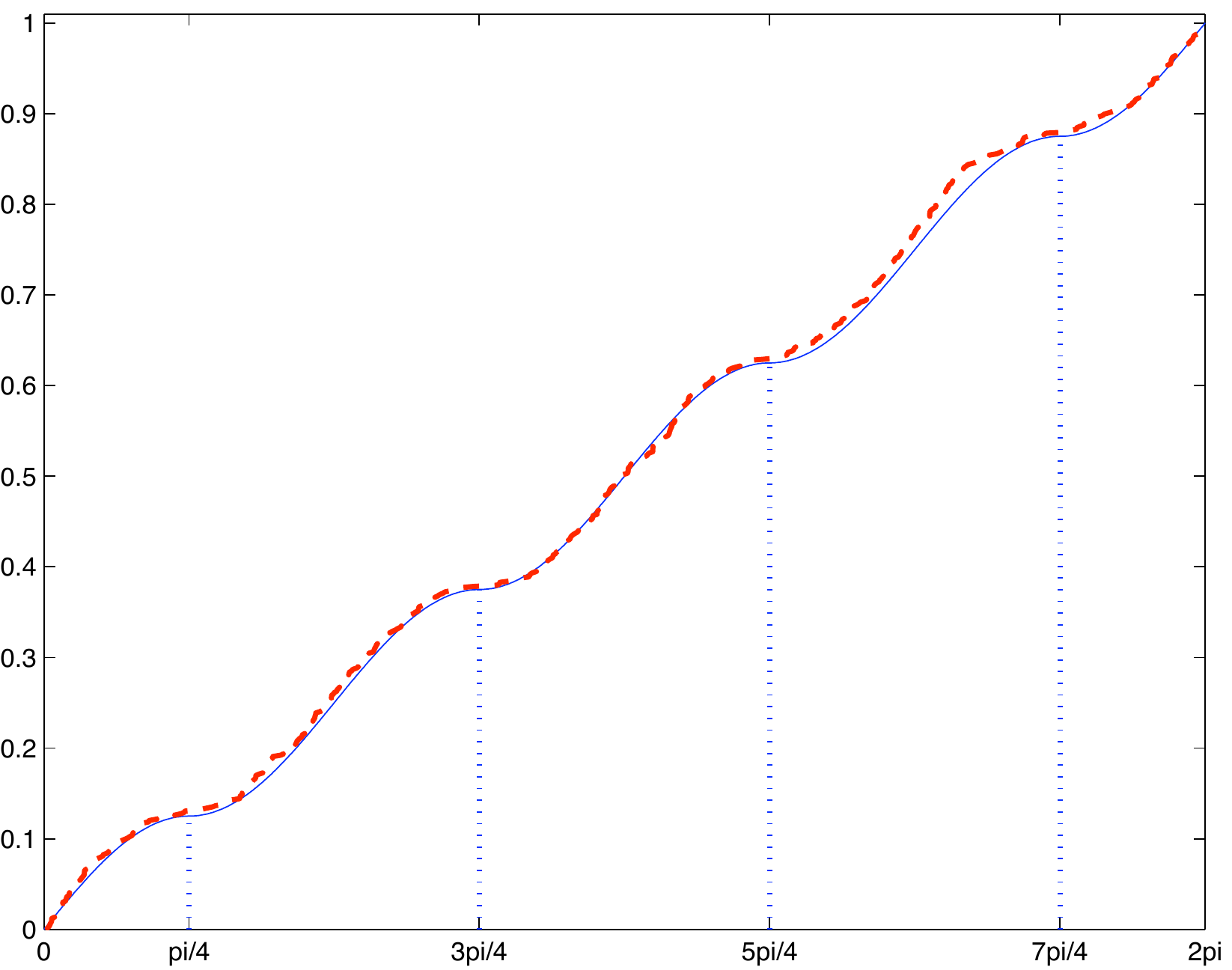} 
\end{center}
\caption{Estimates of the cdf of
spectral measure in Example \ref{ex2}. The solid line is the exact cdf, the
dashed line is the estimated cdf.}
\label{ecdfstable}
\end{figure}

\exemple\label{ex2}
Consider the bivariate strictly stable distribution with $\alpha=0.75$,
$\sigma(S)=1$ and the density of spectral measure
defined by $f(\theta)=\frac{1}{4}|\cos(2\theta)|$, $\theta\in (0,2\pi)$. We
calculated the estimators with $r=0.5$. 
Simulated samples had $50,000$ data vectors. The estimated parameters are
$\hat\alpha=0.74$, $\widehat{\sigma(S)}=0.99$.
The confidence intervals with level $95\%$ are respectively $(\hat\alpha-0.07,
\hat\alpha+0.07)$ and 
$(\widehat{\sigma(S)}-0.11, \widehat{\sigma(S)}+0.13)$. The result of estimation
of the cumulative distribution function (cdf) of
spectral measure is shown in Figure \ref{ecdfstable}.

\vspace{5mm}
The main results of the paper are contained in Section \ref{sect-cons} and
\ref{sect-normasymp}. 
The last section contains the proof of results presented in the previous
sections.

\section{The consistency of estimators}\label{sect-cons}
We assume that $\xi_1, \xi_2,\ldots, \xi_N$ are i.i.d. random $\rd$- valued
vectors with a regularly varying distribution. Our aim is to
estimate the spectral measure $\sigma$ from the sample. The estimators
$\hat\sigma_N(\cdot)$ and $\widehat{\sigma(S)}_N$ 
are defined by (\ref{defestms}) and (\ref{etrmas}).
To establish consistency of these estimators we need two auxiliary results.

Let $X$ be a $\rd$-valued random vector, we denote $G(x)=\reip\{\|X\|>x\}$.
Let $Y_1, Y_2, \ldots$ be the random variables i.i.d. with distribution function
$1-G$ and $Y_{n,1}, Y_{n,2}, \cdots , Y_{n,n}$, $Y_{n,1}\geq Y_{n,2}\geq \cdots
\geq Y_{n,n}$, the corresponding order statistics. Our first result extends a
one-dimensional lemma in \cite{LePage81} (Lemma 1) to $d>1$. 
It says that, for a random variable $X$ satisfying (\ref{vrintro}) with
$\sigma(S)=1$, the vector $b_n^{-1}(Y_{n,1},\ldots,Y_{n,n},0,0,\ldots)$ converge
in distribution to
$(\Gamma_1^{-1/\alpha},\Gamma_2^{-1/\alpha},\ldots)$. The multivariate version
of this lemma is as follows.
\begin{lemme}\label{bny}
If $X\in\mbox{RV}(\alpha,\sigma)$, then the
vector $b_n^{-1}(Y_{n,1},\ldots,Y_{n,n},0,0,\ldots)$ converge in distribution in
$\mathbb{R}^\infty$ to
$\sigma(S)^{1/\alpha}(\Gamma_1^{-1/\alpha},\Gamma_2^{-1/\alpha},\ldots)$, where
$\Gamma_i=\sum^i\limits_{j=1}\lambda_j$ and $\lambda_1, \lambda_2,\ldots$ are
i.i.d. random variables with a standard exponential distribution, i.e. $\reie
(\lambda_{i})=1$.
\end{lemme}
The proof is given in Section \ref{sect-proof}.   

The second result is a variant of the strong law of large numbers for a
triangular array.
\begin{proposition}\label{lemgnc2}
Let $\{X_{m,i}, 1\leq i\leq n\}$ be i.i.d. real random variables for each $m$.
Suppose that the indices $n$ and $m$ satisfy the following relations
\begin{equation}\label{relmnlemgnc2}
n\sim N^r, \;\;\; m\sim N^{1-r} \;\;\; \mbox{as} \;\;\; N\rightarrow\infty
\end{equation}
where $0<r<1$ is a constant and $N\in\mathbb{N}$. If there exists a real number
$k>\frac{2}{r}$ and a constant $M>0$ such that $\reie|X_{m,1}|^{k}\leq
M<\infty$, then
\begin{equation}\label{lgnc2}
\frac{1}{n}\sum_{i=1}^{n}X_{m,i}-\reie
X_{m,1}\xrightarrow[N\rightarrow\infty]{a.s.} 0.
\end{equation}
\end{proposition}

\remarque 
The convergence (\ref{lgnc2}) holds if we replace the condition
$\reie|X_{m,1}|^{k}\leq M$, $k>\frac{2}{r}$ 
by a less restrict  hypothesis $\sum\limits_{m=1}^\infty
m^{-\frac{kr}{2(1-r)}}\reie|X_{m,1}|^k<\infty$.
\vspace{0.5cm}

The proof is given in Section \ref{sect-proof}.
\vspace{5mm}

We consider at first the estimator of the total mass of spectral measure.
\begin{thm}\label{thestmasse}
Let $\xi_1,\ldots,\xi_N$ be i.i.d. $\rd$-valued random vectors such that
the condition (\ref{regulier2}) is satisfied with $b_n=n^{1/\alpha}$. If the
condition (\ref{setting}) holds, then for $0<t<\frac{\alpha r}{2}$,
\begin{equation*}
\widehat{\sigma(S)}_N-\sigma(S)\cvpsn 0,
\end{equation*}
where $\widehat{\sigma(S)}_N$ is defined by (\ref{etrmas}).
\end{thm}
\proof It suffices to prove the following convergence
\begin{equation}\label{cvetrmas}
\frac{1}{n}\sum_{i=1}^n
q_{m,i}^t-\Gamma\left(1-\frac{t}{\alpha}\right)\sigma(S)^{\frac{t}{\alpha}}
\cvpsn 0.
\end{equation}
It follows from Lemma \ref{bny} and the assumption $b_n=n^{1/\alpha}$ that for
all $i$ and $t>0$
\begin{equation}\label{cvqmit}
q_{m,i}^{t}\Rightarrow \frac{\sigma(S)^{t/\alpha}}{\Gamma_1^{t/\alpha}}\;\;
\mbox{as} \ m\rightarrow\infty.
\end{equation}
Since $b_n=n^{1/\alpha}$, the condition (\ref{regulier2}) can be wrote as
\[\lim_{x\rightarrow\infty}x^\alpha\reip\{\|\xi\|>x\}=\sigma(S).\]
Therefore there exist $\delta>0$ and a constant $M>0$ such that for $x>\delta$
we have
\begin{eqnarray*}
\reip\{q_{m,i}^{t}\geq x\}&=&\reip\{M_{m,i}^{(1)}\geq m^{1/\alpha}x^{1/t}\}\\
&=&1-(\reip\{\|\xi\|\leq m^{1/\alpha}x^{1/t}\})^m\\
&\leq&m\reip\{\|\xi\|> m^{1/\alpha}x^{1/t}\}\\
&\leq&Mx^{-\alpha/t}.
\end{eqnarray*}
Taking $x_{0}>\delta$ and $C=x_{0}+Mx_{0}^{1-\alpha/t}\frac{\alpha}{\alpha-t}$
we have for $0<t<\alpha$
\begin{eqnarray}
 \reie q_{m,i}^{t}&=&\int_{\{x<x_{0}\}\cup\{x\geq x_{0}\}}x
d\reip_{q_{m,i}^{t}}(x)\nonumber\\
&\leq&x_{0}+x_{0}\reip\{q_{m,i}^{t}\geq
x_{0}\}+\int_{x_{0}}^\infty\reip\{q_{m,i}^{t}\geq x\}dx\nonumber\\
&\leq&x_{0}+Mx_{0}^{1-\alpha/t}+\int_{x_{0}}^\infty Mx^{-\alpha/t}dx\nonumber\\
&=&C.\label{inqmit}
\end{eqnarray}
Now we choose a real number $0<\delta<1$ and $1<\delta '<\frac{1}{\delta}$. We
set $t=\frac{\alpha r\delta}{2}$ and $k=\frac{2\delta '}{r}$. 
Since $0<tk<\alpha$, it follows from (\ref{inqmit}) that there exists a constant
$C$ such that 
\begin{equation}\label{inqeqmitk}
\reie q_{m,i}^{tk}\leq C
\end{equation}
In combination (\ref{cvqmit}) (\ref{inqeqmitk}) and Proposition \ref{lemgnc2},
we obtain the convergence (\ref{cvetrmas}).
\qed

Let us consider the estimator of the normalized spectral measure defined by
(\ref{defestms}). Note that the random vectors
$\theta_{m,1},\ldots,\theta_{m,n}$ are i.i.d.
in $S$. The following lemma shows the asymptotic property for each
$\theta_{m,i}$, $i=1,\ldots,n$.

\begin{lemme}\label{cvdetheta}
Let $\xi, \xi_1,\ldots,\xi_N$ be i.i.d. $\rd$-valued random vectors and
$\xi\in\mbox{RV}(\alpha,\sigma)$. If $\theta_{m,i}$ is defined by
(\ref{deftheta}), then
\begin{equation*}
\theta_{m,i}\Rightarrow\tilde\sigma \ \mbox{as} \ m\rightarrow\infty
\end{equation*}
for each $i$.
\end{lemme}
\proof For all Borel set $B$ in unite sphere $S$ such that $\sigma(\partial
B)=0$, we have
\begin{eqnarray*}
\reip\{\theta_{m,i}\in B\}&=&\reip\{{\xi_{m,1}}/{\|\xi_{m,1}\|}\in B\}\\
&=&\sum_{k=1}^m\reip\{{\xi_{m,1}}/{\|\xi_{m,1}\|}\in B, \xi_{m,1}={\xi_k}\}\\
&=&m\reip\{{\xi_{m}}/{\|\xi_{m}\|}\in B, \xi_{m,1}=\xi_m\}\\
&=&m\reip\{{\xi_m}/{\|\xi_m\|}\in B, {\|\xi_m\|}\geq{\|\xi_k\|}, \forall
k=1,\ldots,m-1\}\\
&=&m\reip\{{\xi_m}/{\|\xi_m\|}\in B, {\|\xi_m\|}\geq b_m\tau_{m-1}\}\\
&=&\int m\reip\{{\xi_m}/{\|\xi_m\|}\in B, \|\xi_m\|\geq
b_mx\}\reip_{\tau_{m-1}}(dx).
\end{eqnarray*}
where $\reip_{\tau_{m-1}}$ is the distribution of $\tau_{m-1}=\max\limits_{1\leq
k\leq m-1}(\|\xi_k\|b_m^{-1})$.

By (\ref{regulier2}) and Lemma \ref{bny}, the last term converges to
\[\int
\sigma(B)x^{-\alpha}\reip_{\sigma(S)^{1/\alpha}\Gamma_1^{-1/\alpha}}(dx)=\frac{
\sigma(B)}{\sigma(S)}\reie(\Gamma_1)=\tilde\sigma(B).\]
\qed

Therefore for each Borel set $B$ in unit sphere $S$ such that $\sigma(\partial
B)=0$, we have
\[\one_B(\theta_{m,i})\Rightarrow\one_B(\eta) \; \mbox{as} \;
m\rightarrow\infty,\]
where $\eta$ is a random vector with distribution $\tilde\sigma$. This yields
\begin{equation}\label{6.5}
\reie\one_B(\theta_{m,i})\xrightarrow[m\rightarrow\infty]{}\tilde\sigma (B).
\end{equation}

If there exists a constant $r>0$ such that $n\sim N^r$, applying Proposition
\ref{lemgnc2} for the triangular array
$\{\one_B(\theta_{m,1}),\ldots,\one_B(\theta_{m,n})\}$, we have for each fixed
set $B\in\mathcal{B}(S)$ with $\sigma(\partial B)=0$,
\begin{equation}\label{6.6}
\frac{1}{n}\sum_{i=1}^{n}\one_B(\theta_{m,i})-\reie\one_B(\theta_{m,1})\cvpsn 0.
\end{equation}
Together (\ref{6.5}) and (\ref{6.6}) we have the following result.

\begin{thm}\label{connuc}
Let $\xi, \xi_1,\ldots,\xi_N$ be i.i.d. $\rd$-valued random vectors and
$\xi\in \mbox{RV}(\alpha,\sigma)$. If $\hat\sigma_N(\cdot)$ is defined by
(\ref{defestms}) and the condition (\ref{setting}) holds, then $\forall
B\in\mathcal{B}(S)$
with $\sigma(\partial B)=0$,
\begin{equation*}
\hat{\sigma}_N(B)=\frac{1}{n}\sum_{i=1}^{n}\one_B(\theta_{m,i}
)\cvpsn\tilde\sigma(B).
\end{equation*}
\end{thm}

The result is for a fixed set. A stronger convergence can be proved by an
immediate application of the following proposition. 
\begin{proposition}\label{lempscv}
Let $(S,\mathcal{S})$ be a complete separable metric space. Let $\{\sigma_n\}$
be a sequence of random probability measures in $S$. If $\sigma$ is a
probability
measure on $(S,\mathcal{S})$ such that for each set $B\in\mathcal{B}(S)$ with
$\sigma(\partial B)=0$ we have the convergence
$\sigma_n(B)\stackrel{a.s.}{\rightarrow}\sigma(B),$ then
\[\sigma_n\stackrel{a.s.}{\Rightarrow}\sigma \ \mbox{as} \ n\rightarrow\infty.\]
\end{proposition}
The proof is given in Section \ref{sect-proof}.

\begin{corollaire}\label{cvsigmahat}
Under the same assumption of Theorem \ref{connuc}, we have
\[\hat{\sigma}_N\stackrel{a.s.}{\Rightarrow}\tilde\sigma \ \mbox{as} \
N\rightarrow\infty.\]
\end{corollaire}

\section{The asymptotic normality of estimators}\label{sect-normasymp}
In this section we consider the asymptotic normality of the estimators
$\widehat{\sigma(S)}_N$ and $\hat{{\sigma}}_N$. 

\begin{thm}\label{normademasse}
Suppose that random vector $\xi$ satisfies the condition (\ref{cdsdord}) with
$\beta>\alpha+1$, the condition (\ref{setting}) holds. 
If we choose
\begin{eqnarray}\label{condbetamasse}
r&=&{\displaystyle
\frac{3\alpha-4(\beta-1)+\sqrt{16(\beta-1)^2-8\alpha(\beta-1)-7\alpha^2}}{
2\alpha}-\varepsilon,}\;\;\;\mbox{if} \;\;\;
\beta\leq\frac{11}{8}\alpha+1,\nonumber\\
r&=&{\displaystyle \frac{1}{2}-\varepsilon,} \;\;\; \mbox{if} \;\;\;
\beta>\frac{11}{8}\alpha+1,
\end{eqnarray}
where $\varepsilon$ is an arbitrarily small positive constant, then in the
following two cases
\begin{equation}\label{condbetamasse2}
\begin{array}{lll}
a)\;\;  {\displaystyle  0<t<\frac{\alpha r}{4}\wedge 1}&\mbox{if}&
{\displaystyle\alpha+1<\beta\leq \frac{11}{8}\alpha+1 \;\; \mbox{or} \;\;\;
\beta\geq \frac{3}{2}\alpha+1},\\ 
b)\;\;  {\displaystyle 0<t<\frac{3\alpha+2-2\beta}{2}\wedge 1} &\mbox{if}
&{\displaystyle \frac{11}{8}\alpha+1<\beta<\frac{3}{2}\alpha+1},
\end{array}
\end{equation}
we have, as $N\rightarrow\infty$,
\begin{equation}\label{cvnormademasse}
\frac{\displaystyle\sqrt{n}\left(\frac{1}{n}\sum_{i=1}^nq_{m,i}
^t-\Gamma\left(1-\frac{t}{\alpha}\right)\sigma(S)^{t/\alpha}\right)}{
\displaystyle\left(\frac{1}{n}\sum_{i=1}^nq_{m,i}^{2t}-\left(\frac{1}{n}\sum_{
i=1}^nq_{m,i}^t\right)^2\right)^{1/2}}\Rightarrow
\mathcal{N}(0,1).
\end{equation}
\end{thm}

\remarque
Fristedt, \cite{Fristedt72}, proved an asymptotic expansion for the distribution
of the norm of a strictly $\alpha$-stable random vector in $\rd$
\begin{equation}\label{frist}
G(x)=c_1x^{-\alpha}+c_2x^{-2\alpha}+O(x^{-3\alpha}), \ \mbox{as}\;
x\rightarrow\infty.
\end{equation}
That means $\beta=2\alpha$. Therefore the condition of this theorem is satisfied
if $\alpha>1$. If $\alpha>8/5$
the rate of convergence of estimator in $\mathcal{L}_{1}$ is close to $N^{1/4}$.
\vspace{0.5cm}

The proof is given in Section \ref{sect-proof}.
\vspace{0.5cm}

Before considering the asymptotic normality of the estimator of normalized
spectral
measure, we present a {\em strong second-asymptotic relation} : $\forall
B\in\mathcal{B}(S)$ with $\sigma(\partial B)=0$,
\begin{equation}\label{2ndm}
\reip\left\{\frac{\xi}{\|\xi\|}\in B,
\|\xi\|>x\right\}=\sigma(B)x^{-\alpha}+Cx^{-\beta}+o(x^{-\beta}) \;  \mbox{as}
\;
x\rightarrow\infty,
\end{equation}
where $\beta >\alpha >0$. Note that this condition implies the conditions
(\ref{vrintro}) and (\ref{cdsdord}).

Let us denote $\tilde\sigma(B)=b, \  \one_B(\theta_{m,i})=\eta_{m,i},
i=1,2,\ldots,n$.
Then
\[Z_n=\hat{\sigma}_N(B)-\tilde\sigma(B)=\frac{1}{n}\sum_{i=1}^n(\eta_{m,i}
-\tilde\sigma(B)),
\]
\[\sqrt{n}Z_n=\frac{1}{\sqrt{n}}\sum_{i=1}^n(\eta_{m,i}-\tilde\sigma(B))=\frac{1
}{
\sqrt{n}}\sum_{i=1}^n(\eta_{m,i}-\reie\eta_{m,1})+\sqrt{n}(\reie\eta_{m,1}
-\tilde\sigma(B)),\]
\[U_{n}=\frac{1}{\sqrt{n}}\sum_{i=1}^n(\eta_{m,i}-\reie\eta_{m,1}), \  
r_{m}=\reie\eta_{m,1}-\tilde\sigma(B).\]

We set
\[T_N=\displaystyle\frac{\sqrt{n}Z_n}{\left(\displaystyle\frac{1}{n}\sum_{i=1}
^n\eta_{m,i}^2-\left(\displaystyle\frac{1}{n}\sum_{i=1}^n\eta_{m,i}
\right)^2\right)^{1/2}},\]
then the asymptotic property of $\hat{\sigma}_N$ can be described as follows.

\begin{thm}\label{normademc}
Let $\xi, \xi_1,\ldots,\xi_N$ be i.i.d. $\rd$-valued random vectors and
the distribution of $\xi$ satisfies the condition (\ref{2ndm}). If we choose
\[n=N^{2\zeta/(1+2\zeta)-\varepsilon},\;\; m=N^{1/(1+2\zeta)+\varepsilon},\]
where $\zeta=\min(\frac{\beta-\alpha}{\alpha}, 1)$ and $\varepsilon$ is an
arbitrarily small positive constant, then
\begin{equation}\label{29}
T_N\Rightarrow \mathcal{N}(0,1).
\end{equation}
\end{thm}

\remarque 
We can get also asymptotic normality for $\sqrt{n}Z_n$, but the variance of the
limit normal law is $\sigma(B)(1-\sigma(B))$ which we are estimating.

\remarque
If $\xi$ is a strictly $\alpha$-stable random vector in $\rd$, by (\ref{frist})
we get $\beta=2\alpha$, 
thus the asymptotically optimal value of $n$ is approximately
$N^{2/3}$. The rate of convergence of
$\hat{\sigma}_N(B)$ in $\mathcal{L}_{1}$ is close to $N^{1/3}$.
\vspace{0.5cm}

The proof is given in Section \ref{sect-proof}.

\section{Proofs}\label{sect-proof}

\par{\it Preliminary remarks}\\
We recall the definition of a regularly varying function. We say that $L$ is
{\em a regular varying function of index $\alpha $ at infinity (respectively at
origin)} and we denote $ L \in R_\alpha $ (respectively $ L \in R_\alpha (0+) $)
if
\[\frac{L(\lambda x)}{L(x)}\rightarrow x^\alpha,
\ \mbox{as} \  x\rightarrow\infty \; (x\rightarrow 0) \ \mbox{for all} \
\lambda>0.\] 

Let $X$ be $\rd$-valued random vector satisfying the regular variation condition
(\ref{regulier2}). We denote $G(x)=\reip\{\|X\|>x\}$. Then
\begin{equation}\label{nG}
nG(b_nx)\rightarrow \sigma(S)x^{-\alpha},\;  \mbox{as} \; n\rightarrow \infty,
\; \mbox{for all} \; x>0.
\end{equation}
For positive fixed $x$, we choose $n$ the smallest integer such that
$b_{n+1}>x$. Then $b_n\leq x<b_{n+1}$ and for a non-creasing function $G$ we
have
\[\frac{G(\lambda b_{n+1})}{G(b_n)}\leq\frac{G(\lambda
x)}{G(x)}\leq\frac{G(\lambda b_{n})}{G(b_{n+1})},\; \; \mbox{for all}\;
\lambda>0.\]
By (\ref{nG}) we have $nG(b_n)\rightarrow\sigma(S)$, then
\[\frac{G(\lambda x)}{G(x)}\rightarrow\lambda^{-\alpha}, \; \mbox{as} \;
x\rightarrow\infty\;  \mbox{for all}\; \lambda>0.\]
We deduce that $G\in R_{-\alpha}$, which allow us write the following
equivalence
\begin{equation}\label{G}
G(x)\sim x^{-\alpha}L(x),
\end{equation}
where $L(x)$ is a slowly varying function.

We recall a well known result on the asymptotic inverse of a regular varying
function.
\begin{thm}\label{variation} (\cite{Bingham87} Th. 1.5.12)
Let $f\in R_\alpha$ with $\alpha >0$, then $\exists g(x)\in R_{1/\alpha}$ such
that the following relation holds 
\begin{equation}\label{equiva}
f(g(x))\sim g(f(x))\sim x \; \mbox{as} \; x\rightarrow\infty.
\end{equation}
Here $g$ ({\em the asymptotic inverse} of $f$) is defined uniquely up to
asymptotic equivalence, and a version of $g$ is 
\[f^{\leftarrow}(x)=\inf\{y: f(y)\leq x\}.\]
\end{thm}

We denote
\begin{equation}\label{eqvauxthm1}
f(x)=\frac{1}{G(x)}\sim x^\alpha\frac{1}{L(x)},
\end{equation}
then $f(x)\in R_\alpha$ with $\alpha>0$. By applying the previous theorem, we
obtain the inverse $g(x)$ of $f(x)$ in the following form,
\[g(x)=x^{1/\alpha}L^\sharp(x)\]
where the slowly varying function $L^\sharp$ verifies the following relation 
\begin{equation}\label{l1}
L(x)^{-1/\alpha}L^\sharp(f(x))\rightarrow 1,
\end{equation}
and
\begin{equation*}
L(g(x))^{-1}L^\sharp(x)^\alpha\rightarrow 1,\ x\rightarrow\infty .
\end{equation*}
By (\ref{equiva}) and (\ref{eqvauxthm1}) we have
\begin{equation}\label{equideg}
G\left(g\left(\frac{1}{x}\right)\right)\sim x, \ x\rightarrow 0.
\end{equation}
Defining the generalized inverse
\[G^{-1}(x):=\inf\{y: G(y)< x\},\] 
we can prove that
\begin{equation}\label{equideg2}
G(G^{-1}(x))\sim x, \ x\rightarrow 0.
\end{equation}
For this we choose $\lambda>1$, $A>1$, $\delta\in(0,\infty)$, then by the
theorem of Potter (Th. 1.5.6 \cite{Bingham87} page 25) there exists $u_0$ such
that
\[A^{-1}\lambda^{-\alpha-\delta}G(v)\leq G(u)\leq A\lambda^{\alpha+\delta}G(v),
\; \forall v\in[\lambda^{-1}u, \lambda u], \;\; u\geq u_0.\]
We take $x$ small enough such that $G^{-1}(x)\geq u_0$, then by the definition
of $G^{-1}$ there exists $y\in[\lambda^{-1}G^{-1}(x),G^{-1}(x)]$ such that
$G(y)\geq x$, and there exists $y'\in[G^{-1}(x),\lambda G^{-1}(x)]$ such that
$G(y')< x$. Taking $G^{-1}(x)$ for $u$, $y$ and $y'$ for $v$, we get
\[A^{-1}\lambda^{-\alpha -\delta}G(y)\leq G(G^{-1}(x))\leq
A\lambda^{\alpha+\delta}G(y').\]
Hence $\limsup$ and $\liminf$ of $G(G^{-1}(x))/x$ are between
$A\lambda^{\alpha+\delta}$ and $A^{-1}\lambda^{-\alpha-\delta}$  as
$x\rightarrow \infty$. Taking $A$, $\lambda \downarrow 1$, we have
$G(G^{-1}(x))/x\rightarrow 1$. 

The relations (\ref{equideg}) and (\ref{equideg2}) give immediately
\[G^{-1}(x)\sim g(1/x), \ x\rightarrow 0.\]  
Thus we have the equivalent expression of the inverse of $G(x)$:
\begin{equation}\label{G-1}
G^{-1}(x)\sim x^{-1/\alpha}L^\sharp(1/x)\in R_{-1/\alpha}(0+).
\end{equation}

\begin{lemme}
Let $X$ be a $\rd$-valued random vector, $G(x)=\reip\{\|X\|>x\}$. If
$X\in\mbox{RV}(\alpha,\sigma)$, then for each
$i=1,2,\ldots,$
\begin{equation}\label{inverg}  
b_n^{-1}G^{-1}\left(\frac{\Gamma_i}{\Gamma_{n+1}}\right)\cvs\sigma(S)^{1/\alpha}
\Gamma_i^{-1/\alpha} \ \mbox{with probability 1.}
\end{equation}
\end{lemme}

\proof We recall (\ref{nG}) $nG(b_nx)\rightarrow\sigma(S)x^{-\alpha}$ which
implies 
\[G(x_n)\sim\frac{\sigma(S)}{n}\left(\frac{b_n}{x_n}\right)^{\alpha}, \;
n\rightarrow\infty\] 
where $x_n=b_n x$, $x_n\rightarrow\infty$, as $n\rightarrow\infty$. By replacing
the left term in the previous formula by (\ref{G}), we get an equivalent
expression of $b_n$ in terms of $L(x)$:
\begin{equation}\label{bnva}
b_n\sim \left(\frac{nL(x_n)}{\sigma(S)}\right)^{1/\alpha}, \
n\rightarrow\infty.
\end{equation}

Considering (\ref{G-1}), we have an equivalent expression with probability 1 for
each $i$,
\begin{equation}\label{g-1}
G^{-1}\left(\frac{\Gamma_i}{\Gamma_{n+1}}\right)\sim\left(\frac{\Gamma_i}{
\Gamma_{n+1}}\right)^{-1/\alpha}L^\sharp\left(\frac{\Gamma_{n+1}}{\Gamma_{i}}
\right),
\ n\rightarrow\infty
\end{equation}
where $L^\sharp$ satisfies (\ref{l1}) which means
\begin{equation}\label{l2}
L(x_n)^{-1/\alpha}L^\sharp\left(\frac{\Gamma_{n+1}}{\Gamma_{i}}
\right)\rightarrow
1,\ n\rightarrow\infty .
\end{equation}
Collecting (\ref{bnva})-(\ref{l2}) we deduce that with probability $1$
\[b_n^{-1}G^{-1}\left(\frac{\Gamma_i}{\Gamma_{n+1}}\right)\sim\sigma(S)^{
1/\alpha}\left(\frac{\Gamma_{n+1}}{n\Gamma_i}\right)^{1/\alpha}L(x_n)^{-1/\alpha
}L^\sharp\left(\frac{\Gamma_{n+1}}{\Gamma_i}\right)\sim
\left(\frac{\sigma(S)}{\Gamma_i}\right)^{1/\alpha},\]
as $n\rightarrow\infty$, since $\Gamma_{n+1}/n\cvps 1$; the lemma is proved.
\qed

\par{\it Proof of Lemma \ref{bny}:}\quad
It is well known that (see, e.g. \cite{Breiman68} Section 13.6)
\begin{equation}\label{yordre}
(Y_{n,1},\ldots,Y_{n,n})\egenloi
\left(G^{-1}\left(\frac{\Gamma_1}{\Gamma_{n+1}}\right),\ldots,G^{-1}\left(\frac{
\Gamma_n}{\Gamma_{n+1}}\right)\right).
\end{equation}
The lemma follows from (\ref{inverg}) and (\ref{yordre}).\qed

\par{\it Proof of Proposition \ref{lemgnc2}:}\quad
Denote $Y_{m,i}=X_{m,i}-\reie X_{m,1}$, then the random variables $\{Y_{m,i},
1\leq i\leq n\}$ are centered and i.i.d.. We have
\[\reie|Y_{m,1}|^{k}\leq \reie (|X_{m,1}|+|\reie
X_{m,1}|)^{k}\leq\reie(2^{k-1}(|X_{m,1}|^{k}+|\reie X_{m,1}|^{k}))\leq
2^{k}\reie|X_{m,1}|^{k}.\]
It is well  known (see \cite{Rosenthal70}) that for $k\geq 2$ we have
\[\reie\left|\sum_{i=1}^nY_{m,i}\right|^{k}\leq
c(k)n^{k/2}\reie|Y_{m,1}|^{k},\] 
where $c(k)$ is a positive constant depending only on $k$. It follows from the
condition (\ref{relmnlemgnc2}) that there exists a constant $C>0$ such that
$n\geq CN^r$. Hence for all $\varepsilon >0$ we have
\[\reip\left\{\left|\frac{1}{n}\sum_{i=1}^nY_{m,i}\right|>\varepsilon\right\}
\leq\frac{\reie\left|{\displaystyle
\sum_{i=1}^nY_{m,i}}\right|^{k}}{n^{k}\varepsilon^{k}}\leq
\frac{2^{k}c(k)\reie|X_{m,1}|^{k}}{C^{\frac{k}{2}}N^{\frac{kr}{2}}\varepsilon^{k
}}=\frac{c_0}{N^{\frac{kr}{2}}\varepsilon^{k}},\]
where $c_0=2^{k}c(k)M/C^{\frac{k}{2}}$. Since $\frac{kr}{2}>1$, we can find a
small enough positive number $\varepsilon'$ such that $\frac{kr}{2}-\varepsilon
'>1$. Taking $\varepsilon=\varepsilon_N=N^{-\frac{\varepsilon '}{k}}$ and
applying the Borel-Cantelli lemma, we have that with probability $1$ and for $N$
large enough
\[\left|\frac{1}{n}\sum_{i=1}^nX_{m,i}-\reie X_{m,1}\right|\leq
N^{-\frac{\varepsilon '}{k}},\]
the proposition is proved. 
\qed

\par{\it Proof of Proposition \ref{lempscv}:}\quad
We denote the collection of all $\sigma$-continuity sets by
\[\mathcal{D}_\sigma=\{B~|~B\in\mathcal{B}(S),\;\sigma(\partial B)=0\}.\] 
Since space $S$ is separable, there exists a countable dense set in $S$, denoted
by
\[W=\{x_1,x_2,\ldots\},\ x_i\in S, \ i=1,2,\ldots.\]
We denote the open ball with centre $x_i$ in $W$ and radius $r$ by 
\[V(x_i,r)=\{x~|~ x\in S, \|x-x_i\|<r\}.\]
Since for each $x_i\in W$ the boundaries
$\partial\{V(x_i,r)\}\subset\{x~|~\|x-x_i\|=r\}$ are disjoints for different
$r$, at most a countable number of them can have positive $\sigma$-measure.
Therefore, there exists a sequence of positive numbers $r_k^i \downarrow 0$ as
$k\rightarrow\infty$ for each $x_i$ such that
\[\mathbb{L}_{i}=\{V(x_i,r_k^i), k=1,2,\ldots\}\subset\mathcal{D}_\sigma.\]
The collection $\mathbb{L}=\bigcup\limits_{x_i\in W}\mathbb{L}_{i}$ is
countable. It is clear that for each $x_i\in W$, the collection $\mathbb{L}_{i}$
is a local base at point $x_i$ for la topology $\mathcal{S}$. Since $W$ is dense
in $S$, $\mathbb{L}$ is a base of $\mathcal{S}$. The $\sigma$-algebra generated
by $\mathbb{L}$, denoted by $\sigma(\mathbb{L})$, is the Borel-field
$\mathcal{B}(S)$.

Now we expand $\mathbb{L}$ by adding the finite intersections of members of
$\mathbb{L}$, we denote
\[\mathcal{L}=\mathbb{L}\cup\left.\left\{\bigcap_{i\in
I}V_{i}~\right|~V_{i}\in\mathbb{L}, I\subset\mathbb{N},
\;\mbox{card}(I)<\infty\right\}.\]
It is clear that $\mathcal{L}$ is still countable and
$\sigma(\mathcal{L})=\mathcal{B}(S)$, moreover
$\mathcal{L}\subset\mathcal{D}_\sigma$. Since
$\sigma_n(B)\stackrel{a.s.}{\rightarrow}\sigma(B)$ for all $B\in\mathcal{B}(S)$
and $\sigma(\partial B)=0$, then $\forall V\in\mathcal{L}$,
$\exists\Lambda_V\subset\Omega$ and $\reip(\Lambda_V)=0$, such that $\forall
\omega\in\Lambda_V^\complement$ we have
\begin{equation}\label{cvsigma}
\sigma_n(\omega,V)\rightarrow\sigma(V).
\end{equation}
If we denote $\Lambda=\bigcup\limits_{V\in\mathcal{L}}\Lambda_V$, then
$\reip(\Lambda)=0$. Moreover $\forall \omega\in\Lambda^\complement$ we have
always the convergence (\ref{cvsigma}) for all $V\in\mathcal{L}$. The collection
$\mathcal{L}$ is closed under the operation of finite intersection. By Theorem
2.2 in \cite{Billingsley68} (page 14) we have $\sigma_n\Rightarrow\sigma$,
$\forall \omega\in\Lambda^\complement$, which implies
$\sigma_n\stackrel{a.s.}{\Rightarrow}\sigma .$ \qed

\par{\it Proof of Theorem \ref{normademasse}:}\quad
We set $X_{m,i}=q_{m,i}^t$, $\mu_{m}=\reie q_{m,i}^t$,
$\sigma_m^2=\mbox{Var}(X_{m,i})$ and
$r_m=\mu_m-\Gamma(1-t/\alpha)\sigma(S)^{t/\alpha}$. 
Therefore the convergence (\ref{cvnormademasse}) holds if we have the following
three relations :
\begin{equation}\label{cvnormademasse1}
\frac{1}{\sqrt{n}}\sum_{i=1}^n(X_{m,i}-\mu_{m})\Rightarrow\mathcal{N}(0,
\sigma^2)
\end{equation}

\begin{equation}\label{cvnormademasse3}
\frac{1}{n}\sum_{i=1}^nX_{m,i}^{2}-\left(\frac{1}{n}\sum_{i=1}^nX_{m,i}
\right)^2\xrightarrow[N\rightarrow\infty]{P}\sigma^2.
\end{equation}
and
\begin{equation}\label{cvnormademasse2}
\sqrt{n}r_m\rightarrow 0
\end{equation}

By the similar method to (\ref{inqeqmitk}) we can prove the moments $\{\reie
X_{m,i}^4\}$ are uniformly bounded. Hence we have the following convergence
\begin{equation}\label{sigmanormademasse}
\sigma_m^2:=\reie X_{m,i}^{2}-(\reie
X_{m,i})^2\rightarrow\sigma^2=\sigma(S)^{2t/\alpha}
(\Gamma(1-2t/\alpha)-(\Gamma(1-t/\alpha))^2),
\end{equation}
and the Lindeberg's condition, i.e. for all $\varepsilon>0$
\begin{eqnarray*}
&&\lim_{n\rightarrow\infty}\frac{1}{n\sigma_m^2}\sum_{i=1}^n\int_{\{|X_{m,i}
-\mu_m|>\varepsilon\sqrt{n}\sigma_m\}}(X_{m,i}-\mu_m)^2d\reip\\
&=&\lim_{n\rightarrow\infty}\frac{1}{\sigma_m^2}\int_{\{|X_{m,1}
-\mu_m|>\varepsilon\sqrt{n}\sigma_m\}}(X_{m,1}-\mu_m)^2d\reip\\
&\leq&\lim_{n\rightarrow\infty}\frac{1}{\sigma_m^2}(\reie(X_{m,1}-\mu_m)^4)^{1/2
}(\reip\{|X_{m,1}-\mu_m|>\varepsilon\sqrt{n}\sigma_m\})^{1/2}\\
&=&0.
\end{eqnarray*}
The convergence (\ref{cvnormademasse1}) follows from the central limit theorem
applied to triangular array  
$\{X_{m,i}, 1\leq i\leq n\}$. By Proposition \ref{lemgnc2} and inequality
(\ref{inqmit}), if $0<t<\frac{\alpha r}{4}$ we have the following convergences 
\begin{equation*}
\frac{1}{n}\sum_{i=1}^nX_{m,i}^{2}-\reie X_{m,i}^{2}\cvpsn 0
\end{equation*}
and
\begin{equation*}
\frac{1}{n}\sum_{i=1}^nX_{m,i}-\mu_m\cvpsn 0.
\end{equation*}
It follows
\[\frac{1}{n}\sum_{i=1}^n
X_{m,i}^{2}-\left(\frac{1}{n}\sum_{i=1}^nX_{m,i}\right)^2-(\reie
X_{m,i}^{2}-\mu_m^2)\cvpsn 0,\]
considering (\ref{sigmanormademasse}) we obtain (\ref{cvnormademasse3}). Now it
remains to verify(\ref{cvnormademasse2}).

\begin{lemme}\label{rmmasse}
If $\xi$ satisfies the condition (\ref{cdsdord}) with $\beta>\alpha+1$ and
$0<t<1$, then
\begin{equation}\label{11}
|r_m|\leq Cm^{-\zeta}
\end{equation}
where $\zeta=\min(\frac{1}{2},\frac{t+\beta-\alpha-1}{\alpha})$ and $C$
depending only on $c_1, c_2,
\alpha$, $\beta$ and $t$.
\end{lemme}

\proof By relation (\ref{cdsdord}) if $x$ is sufficiently large we have
\[G(x)=c_1x^{-\alpha}+c_2x^{-\beta}+o(x^{-\beta}),\]
where $c_{1}=\sigma(S)$. It is possible to write the inverse function for small
value of  $t$,
\begin{equation*}
G^{-1}(t)=\sigma(S)^{1/\alpha}t^{-1/\alpha}+bt^{s}+O(t^{s+(\beta-\alpha)/\alpha}
)
\end{equation*}
with $b=\alpha^{-1}c_2\sigma(S)^{(1-\beta)/\alpha}$ et
$s=(\beta-\alpha-1)/\alpha$. It follows that for small $\delta>0$ and
$0<t<\delta$
\[G^{-1}(t)-\sigma(S)^{1/\alpha}t^{-1/\alpha}=t^s(b+O(t^{(\beta-\alpha)/\alpha}
)).\]
We choose $\delta$ such that $|O(t^{(\beta-\alpha)/\alpha})|\leq |b|$ (this
gives us $\delta=o(|b|^{\frac{\alpha}{\beta-\alpha}})$), then we can write
\begin{equation}\label{ingg}
|G^{-1}(t)-\sigma(S)^{1/\alpha}t^{-1/\alpha}|\leq 2|b|t^s.
\end{equation}

Let
\[R^{m+1}_+=\{\bar{x}=(x_1,\ldots,x_{m+1}) : x_i\geq 0,i=1,\ldots,m+1\},\]
\[\Sigma_m=x_1+\cdots +x_m,\]
\[A=\left\{\bar{x}\in R^{m+1}_+ : \frac{x_1}{\Sigma_{m+1}}\geq\delta\right\},
A^\complement=R_+^{m+1}\backslash A,\]
where $\delta$ is chosen for that (\ref{ingg}) holds. By its definition and the
relation (\ref{yordre}) the random variable 
$q_{m,i}=M_{m,i}^{(1)}/m^{1/\alpha}$ have the distribution 
\[G^{-1}\left(\frac{\Gamma_1}{\Gamma_{m+1}}\right)\Big/m^{1/\alpha}\]
then
\begin{equation}\label{14}
|r_m|=\left|\reie\left(G^{-1}\left(\frac{\Gamma_1}{\Gamma_{m+1}}
\right)\right)^t\Big/m^{t/\alpha}-\reie\left(\frac{\sigma(S)}{\Gamma_{1}}
\right)^{t/\alpha}\right|\leq
I_1+I_2+I_3+I_{4}, 
\end{equation}
where
\begin{eqnarray*}
I_1&=&\int_{A}m^{-t/\alpha}\left|\left(G^{-1}\left(\frac{x_1}{\Sigma_{m+1}}
\right)\right)^t\right|\exp(-\Sigma_{m+1})d\bar{x},\\
I_2&=&\int_{A^\complement}m^{-t/\alpha}\left|\left(G^{-1}\left(\frac{x_1}{
\Sigma_{m+1}}\right)\right)^t-\left(\frac{x_1}{\Sigma_{m+1}}\right)^{-t/\alpha}
\sigma(S)^{t/\alpha}\right|\exp(-\Sigma_{m+1})d\bar{x},\\
I_3&=&\int_{A^\complement}m^{-t/\alpha}\sigma(S)^{t/\alpha}x_{1}^{-t/\alpha}
\left|\Sigma_{m+1}^{t/\alpha}-m^{t/\alpha}\right|\exp(-\Sigma_{m+1})d\bar{x},\\
\mbox{and}&&\\
I_4&=&\int_{A}\sigma(S)^{t/\alpha}x_{1}^{-t/\alpha}\exp(-\Sigma_{m+1})d\bar{x}.
\end{eqnarray*}
Since in $A$ we have $G^{-1}(x_{1}/\Sigma_{m+1})\leq G^{-1}(\delta)$, for all
$\tau_{1}>0$ (which will be chosen later)
\begin{equation}\label{ingi1}
I_1\leq
m^{-t/\alpha}(G^{-1}(\delta))^t\reip\left\{\frac{\Gamma_1}{\Gamma_{m+1}}
\geq\delta\right\}\leq
m^{-t/\alpha}C_{1}\reie\left(\frac{\Gamma_1}{\Gamma_{m+1}}\right)^{\tau_{1}}
\end{equation}
where $C_{1}=(G^{-1}(\delta))^t\delta^{-\tau_{1}}$. In a similar way we estimate
$I_4$, for all $\tau_{2}>0$,
\begin{equation}\label{ingi4}
I_4\leq
\sigma(S)^{t/\alpha}(\reie\Gamma_{1}^{-2t/\alpha})^{1/2}\left(\reip\left\{\frac{
\Gamma_1}{\Gamma_{m+1}}\geq\delta\right\}\right)^{1/2}\leq
C_{2}\reie\left(\frac{\Gamma_1}{\Gamma_{m+1}}\right)^{\tau_{2}}
\end{equation}
where $C_{2}=\sigma(S)^{t/\alpha}(\Gamma(1-2t/\alpha))^{1/2}\delta^{-\tau_{2}}$.
The function $f(x)=x^t, 0<t<1$ is Lipschitz continuous on $[0,\delta)$.
Considering (\ref{ingg}) we have
\begin{eqnarray}
I_{2}&\leq&\int_{A^\complement}m^{-t/\alpha}C_{3}\left|G^{-1}\left(\frac{x_1}{
\Sigma_{m+1}}\right)-\left(\frac{x_1}{\Sigma_{m+1}}\right)^{-1/\alpha}\sigma(S)^
{1/\alpha}\right|\exp(-\Sigma_{m+1})d\bar{x}\nonumber\\
&\leq&2|b|C_{3}m^{-t/\alpha}\int_{A^\complement}\left(\frac{x_{1}}{\Sigma_{m+1}}
\right)^s\exp(-\Sigma_{m+1})d\bar{x}\nonumber\\
&\leq&2|b|C_{3}m^{-t/\alpha}\reie\left(\frac{\Gamma_{1}}{\Gamma_{m+1}}
\right)^s\label{ingi2}
\end{eqnarray}
where
$C_{3}=t(\min(G^{-1}(\delta),\delta^{-1/\alpha}\sigma(S)^{1/\alpha}))^{t-1}$. 
It is well known (see for example \cite{Breiman68}) that the m-dimensional
random vector 
$\displaystyle{\left(\frac{\Gamma_{1}}{\Gamma_{m+1}},\ldots,\frac{\Gamma_{m}}{
\Gamma_{m+1}}\right)}$ 
has the same density as the order statistics vector of random variables
uniformly distributed on  $[0,1)$. 
In particular the random variable
$\displaystyle{\frac{\Gamma_{1}}{\Gamma_{m+1}}}$ has the following density
(suppose that $m\geq 2$)
\[g(x)=\left\{\begin{array}{ll} m(1-x)^{m-1},&\mbox{if} \ 0\leq x\leq
1,\\0,&\mbox{otherwise}.\end{array}\right.\]
Thus
\begin{equation}\label{inggammas}
\reie\left(\frac{\Gamma_1}{\Gamma_{m+1}}\right)^s=m\int_0^1x^s(1-x)^{m-1}
dx=mB(s+1,m)\leq 4\Gamma(s+1)m^{-s}.
\end{equation}
By Cauchy Schwarz's inequality we have 
\begin{equation}\label{ingi3}
I_{3}\leq
m^{-t/\alpha}\sigma(S)^{t/\alpha}(\reie\Gamma_1^{-2t/\alpha})^{1/2}
(\reie(\Gamma_{m}^{t/\alpha}-m^{t/\alpha})^2)^{1/2}.
\end{equation}
It remains to evaluate 
\begin{eqnarray*}
\reie(\Gamma_{m}^{t/\alpha}-m^{t/\alpha})^2&=&\frac{\Gamma(m+1+2t/\alpha)}{
\Gamma(m+1)}-2m^{t/\alpha}\frac{\Gamma(m+1+t/\alpha)}{\Gamma(m+1)}+m^{2t/\alpha}
.
\end{eqnarray*}
Using the approximation of Gamma-function 
\[\Gamma(x)=\sqrt{2\pi}x^{x-1/2}e^{-x}\left(1+\frac{1}{12x}+\frac{1}{288x^2}
+o\left(\frac{1}{x^2}\right)\right), \;\;\; x\rightarrow\infty,\]
we obtain 
\begin{equation}\label{egespi3}
\reie(\Gamma_{m}^{t/\alpha}-m^{t/\alpha})^2=A_m(m+1+2t/\alpha)^{2t/\alpha}
H_1-2B_mm^{t/\alpha}(m+1+t/\alpha)^{t/\alpha}H_2+m^{2t/\alpha},
\end{equation}
where
\[A_m=\left(1+\frac{2t/\alpha}{m+1}\right)^{m+1}\left(1+\frac{2t/\alpha}{m+1}
\right)^{-1/2}e^{-2t/\alpha}, \]
\[B_m=\left(1+\frac{t/\alpha}{m+1}\right)^{m+1}\left(1+\frac{t/\alpha}{m+1}
\right)^{-1/2}e^{-t/\alpha},\]
\begin{equation}\label{devh}
H_1=1-\frac{t}{6\alpha m^2}+o\left(\frac{1}{m^2}\right)
\;\;\;\mbox{et}\;\;\;H_2=1-\frac{t}{12\alpha m^2}+o\left(\frac{1}{m^2}\right).
\end{equation}
By a simple calculation we obtain the following expansion
\begin{equation}\label{devam}
A_m=1-\left(\frac{t}{\alpha}+\frac{2t^2}{\alpha^2}\right)\frac{1}{m+1}
+o\left(\frac{1}{m}\right),
\end{equation}
\begin{equation}\label{devbm}
B_m=1-\left(\frac{t}{2\alpha}+\frac{t^2}{2\alpha^2}\right)\frac{1}{m+1}
+o\left(\frac{1}{m}\right).
\end{equation}
Considering (\ref{egespi3})-(\ref{devbm}) we have
\[\reie(\Gamma_{m}^{t/\alpha}-m^{t/\alpha})^2=\frac{t^2}{\alpha^2}m^{2t/\alpha-1
}+o(m^{2t/\alpha-1}).\]
Therefore it follows from (\ref{ingi3}) that there exists a positive constant
$C_4$ such that
\begin{equation}\label{ingi32}
I_3\leq C_4m^{-1/2}.
\end{equation}
Collecting estimates (\ref{14})-(\ref{inggammas}), (\ref{ingi32}) and choosing
$\tau_{1} =(\beta-\alpha-1)/\alpha$ and $\tau_{2} =(t+\beta-\alpha-1)/\alpha$ 
in (\ref{ingi1}) and (\ref{ingi4}), we obtain the estimate (\ref{11}) with
$C=C(\alpha,\beta,c_1,c_2, t)$. The lemma is proved.
 \qed

Having the estimate for $r_m$, and taking $n\sim N^{r}$ with
$r<\frac{2\zeta}{1+2\zeta}$ we obtain $\sqrt{n}r_m\rightarrow 0$. 
By a simple calculation we have the relations (\ref{condbetamasse}) and
(\ref{condbetamasse2}).
\qed

\par{\it Proof of Theorem \ref{normademc}:}\quad
We denote $c^2=b(1-b)$. In order to prove (\ref{29}) it suffices to show the
following relations,
\begin{equation}\label{30}
U_n\cvdn \mathcal{N}(0,c^2),
\end{equation}

\begin{equation}\label{31}
\sqrt{n}r_m\rightarrow 0,
\end{equation}

\begin{equation}\label{32}
\displaystyle\frac{1}{n}\sum_{i=1}^n\eta_{m,i}^2-\left(\displaystyle\frac{1}{n}
\sum_{i=1}^n\eta_{m,i}\right)^2\xrightarrow[N\rightarrow\infty]{P}
c^2.
\end{equation}

Since $0<\one_B(\theta_{m,i})\leq 1$, the moments
$\reie|\one_B(\theta_{m,i})-\reie \one_B(\theta_{m,i})|^k$ are uniformly bounded
for all $k>0$. The limit variance is 
\[c_m^2:=\reie(\eta_{m,1}-\reie\eta_{m,1})^2=\reie\eta^2_{m,1}-(\reie\eta_{m,1}
)^2=b+r_m-(b+r_m)^2\rightarrow c^2,\]
if we have (\ref{31}). We consider the Lindeberg's condition: for all
$\varepsilon>0$,
\begin{eqnarray}\label{ingkappa}
&&\lim_{n\rightarrow\infty}\frac{1}{nc_m^2}\sum_{i=1}^n\int_{\{|\one_B(\theta_{m
,i})-\reie
\one_B(\theta_{m,i})|>\varepsilon\sqrt{n}c_m\}}(\one_B(\theta_{m,i})-\reie
\one_B(\theta_{m,i}))^2d\reip\nonumber\\
&=&\lim_{n\rightarrow\infty}\frac{1}{c_m^2}\int_{\{|\one_B(\theta_{m,1})-\reie
\one_B(\theta_{m,1})|>\varepsilon\sqrt{n}c_m\}}(\one_B(\theta_{m,1})-\reie
\one_B(\theta_{m,1}))^2d\reip\nonumber\\
&\leq&\lim_{n\rightarrow\infty}\frac{1}{c_m^2}(\reie(\one_B(\theta_{m,1})-\reie
\one_B(\theta_{m,1}))^4)^{1/2}(\reip\{|\eta_{m,1}-\reie\eta_{m,1}
|>\varepsilon\sqrt{n}
c_m\})^{1/2}\nonumber\\
&=&0.
\end{eqnarray}
The relation (\ref{30}) follows from the application of central limit theorem
for the triangular array $\{\one_B(\theta_{m,i})-\reie \one_B(\theta_{m,i})\}$.

It is easy to see that (\ref{32}) follows from the application of Proposition
\ref{lemgnc2} to triangular array $\{\eta_{m,i}^2\}$. It remains to establish
the relation (\ref{31}). In fact, we have $|r_m|\leq C\max(m^{-1},
m^{-(\beta-\alpha)/\alpha})$ from the following lemma. Taking an arbitrarily
small positive constant $\varepsilon$ and
\[n=N^{\frac{2\zeta}{1+2\zeta}-\varepsilon}, \;
m=N^{\frac{1}{1+2\zeta}+\varepsilon},\]
we get $\sqrt{n}r_{m}\rightarrow 0$, the theorem is proved.\qed

\begin{lemme}\label{1}
If the condition (\ref{2ndm}) is satisfied, then
\begin{equation}\label{33}
|r_m|\leq C\max(m^{-1}, m^{-(\beta-\alpha)/\alpha}).
\end{equation}
\end{lemme}
\proof In order to prove (\ref{33}), we need to show that
\begin{equation}\label{34}
\reip\{\theta_{m,i}\in B\}=\sigma(B)+R_m,
\end{equation}
with the remainder term $R_m=O(\max(m^{-1}, m^{-(\beta-\alpha)/\alpha}))$. Let
us denote
\[G_m(x)=\reip\left\{\max_{1\leq i\leq m-1}\|\xi_i\|\leq x\right\}.\]
Using the definition of $\theta_{m,i}$, it is not difficult to see that
\[\reip\{\theta_{m,i}\in
B\}=m\int_0^\infty\reip\left\{\frac{\xi_1}{\|\xi_1\|}\in B,
\|\xi_1\|>r\right\}G_m(dr).\]
Let $\widetilde{G}_m(x)=G_m(xm^{1/\alpha}).$ Assumption (\ref{2ndm}) implies (we
suppose that $\sigma(S)=1$) that for large $s$,
\begin{equation*}\label{35}
\reip\{\|\xi\|>s\}=s^{-\alpha}+Cs^{-\beta}+o(s^{-\beta}).
\end{equation*}
Therefore, it is easy to get the relation
\[\lim_{m\rightarrow\infty}\widetilde{G}_m(x)=G_0(x)=\left\{\begin{array}{ll}
\exp({-x^{-\alpha}}),&
x>0,\\
0,& x\leq 0.\end{array}\right.\] 
Using (\ref{2ndm}) and the fact that $\int_0^\infty y^{-\alpha}dG_0(y)=1$, we
have (\ref{34}) with $R_m=\sum_{i=1}^4 R_{m,i}$, where
\begin{eqnarray*}
R_{m,1}&=&m\int_0^s\reip\left\{\frac{\xi_1}{\|\xi_1\|}\in B,
\|\xi_1\|>r\right\}dG_m(r),\\
R_{m,2}&=&-\sigma(B)\int_o^{s'}y^{-\alpha}dG_0(y),\\
R_{m,3}&=&\sigma(B)\int_{s'}^\infty
y^{-\alpha}d(\widetilde{G}_m(y)-G_0(y)),\\
R_{m,4}&=&Cm^{-(\beta-\alpha)/\alpha}\int_{s'}^\infty
y^{-\beta}d\widetilde{G}_m(y).
\end{eqnarray*}
Here $s'=sm^{-1/\alpha}$ and we shall choose $s$ later. It is easy to see that
\[R_{m,1}\leq m(1-\reip\{\|\xi_1\|>s\})^m=m(1-hm^{-1})^m\leq
me^{-\frac{1}{2}h},\]
where $h=h(m,s)=m\reip\{\|\xi_1\|>s\}\geq\frac{1}{2}ms^{-\alpha}$. We have used
(\ref{35}) for the last inequality. Thus, if we choose
\[s=\left(\frac{m}{K\ln m}\right)^{1/\alpha},\]
with sufficiently large $K$, then we get
\begin{equation}\label{36}
R_{m,1}=o(m^{-1}).
\end{equation} 
Simple calculations show that
\begin{equation}\label{37}
R_{m,2}=o(m^{-1}).
\end{equation}
The main remainder term is $R_{m,3}$ and to estimate it we must first estimate
the difference $\widetilde{G}_m(y)-G_0(y)$. A rather simple expansion of
logarithmic function gives the following estimates which are sufficient for our
purposes.

\begin{lemme}\label{2} (\cite{Davydov99} Lemma 2)
Let $\xi_i, i\geq 1$ be i.i.d. random vectors satisfying (\ref{35}). Then for
$y>cm^{-1/\alpha}$
\begin{equation}\label{38}|\widetilde{G}_m(y)-G_0(y)|\leq
C(\alpha,\beta)\exp(-y^{-\alpha})(m^{-(\beta-\alpha)/\alpha}y^{-\beta}+m^{-1}y^{
-2\alpha})
\end{equation}
and
\begin{equation*}
\sup_y|\widetilde{G}_m(y)-G_0(y)|=C(\alpha,\beta)\max(m^{-1},m^{
-(\beta-\alpha)/\alpha}).
\end{equation*}
\end{lemme}

Now we can estimate the term $R_{m,3}$. Integrating by parts, we get
\begin{equation}\label{40}
|R_{m,3}|=\sigma (B)(R_{m,3}^{(1)}+R_{m,3}^{(2)}),
\end{equation}
where
\[R_{m,3}^{(1)}=s'^{-\alpha}|\widetilde{G}_m(s')-G_0(s')|,\]
\[R_{m,3}^{(2)}=\alpha\int_{s'}^\infty|\widetilde{G}_m(y)-G_0(y)|y^{-\alpha-1}
dy.\]
Since $s'=(K\ln m)^{-1/\alpha}>C m^{-1/\alpha}$, we can use (\ref{38}) to
estimate both quantities $R_{m,3}^{(i)}, i=1,2$. After some simple calculations,
we get
\[R_{m,3}^{(1)}=o(m^{-1}),\]
\[R_{m,3}^{(2)}\leq C(\alpha,\beta)\max(m^{-1},m^{-(\beta-\alpha)/\alpha}).\]

In a similar way we estimate $R_{m,4}$:
\begin{equation}\label{42}
R_{m,4}=Cm^{-(\beta-\alpha)/\alpha}\int_{s'}^\infty
y^{-\beta}d\widetilde{G}_m(y)=Cm^{-(\beta-\alpha)/\alpha}(R_{m,4}^{(1)}+R_{m,4}^
{(2)}),
\end{equation}
where
\[R_{m,4}^{(1)}=\int_{s'}^\infty y^{-\beta}dG_0(y),\]
\[R_{m,4}^{(2)}=\int_s^\infty y^{-\beta}d(\widetilde{G}_m(y)-G_0(y)).\]
It is easy to see that
\begin{equation}\label{43}
R_{m,4}^{(1)}\leq C(\alpha,\beta)
\end{equation}
and $R_{m,4}^{(2)}$ can be estimated in a similar way to $R_{m,3}$:
\begin{equation}\label{44}
R_{m,4}^{(2)}\leq C(\alpha,\beta)\max(m^{-1},
m^{-(\beta-\alpha)/\alpha}).
\end{equation}
Collecting (\ref{36}), (\ref{37}), and (\ref{40})-(\ref{44}) we get (\ref{34}). 
\qed

\bibliographystyle{plain}

\end{document}